\documentclass[11pt,reqno]{amsart}
\usepackage{amsfonts}
\usepackage{amsmath}
\usepackage{amsthm}
 
\usepackage{hyperref}
\usepackage[capitalise,nameinlink,noabbrev]{cleveref}

\usepackage{amssymb}
\usepackage[all]{xy}

\DeclareFontFamily{OT1}{pzc}{}
\DeclareFontShape{OT1}{pzc}{m}{it}{<-> s * [1.200] pzcmi7t}{}
\DeclareMathAlphabet{\mathpzc}{OT1}{pzc}{m}{it}

\newcommand{\mathcat}{\mathpzc}

\def\Out{\mathrm{Out}} 
\def\Aut{\mathrm{Aut}}

\def\Id{\mathrm{Id}}

\def\Ho{\mathrm H}
\def\Id{\mathrm{Id}}
\def\kappaQ{\kappa}

\def\res{\mathrm{res}}
\def\fiel{\mathfrak k}

\title{Crossed modules}

\author{Johannes Huebschmann  }
\address{
\newline
Universit\'e de Lille - Sciences et Technologies 
\\
D\'epartement de Math\'ematiques\\
\newline CNRS-UMR 8524,
Labex CEMPI (ANR-11-LABX-0007-01)
\\
\newline
59655 Villeneuve d'Ascq Cedex, France\\
\newline
Johannes.Huebschmann@univ-lille.fr
 }

\begin{document}
\thispagestyle{empty}
\begin{abstract}
\noindent
This is an overview of the idea of a crossed module.
For a group, the triple that consists of the group, its group of 
automorphisms,
and  the canonical homomorphism from the group to its group 
of automorphisms constitutes a crossed module.
Crossed modules arise from the identities among the relations
of the presentation of a group, from the extension problem for groups
and, more generally, in low dimensional topology.
Also, the (successful) attempt to extend the idea of a normal extension
of commutative fields to the realm of non-commutative algebras
 leads to crossed modules. Crossed modules appear implicitly in a forgotten
paper by A. Turing which in principle 
settles the extension problem for groups.
Crossed modules make perfect sense for Lie algebras.
\end{abstract}
\maketitle

\centerline
{Dedicated to the memory of Donald J. Collins}

\section{Introduction}
Modern chatbot software poses a threat to the health of our field.
A scholarly article had better pass
through the heads of at least two parties, 
cf. \cite[p.~47]{huff2010lie}.\footnote{\lq\lq 
It appears that the reporter has passed along some words without inquiring
what they mean, and you are expected to read them just 
as uncritically for the happy illusion they give you of having learned 
something.
It is all too reminiscent
of an old definition of the lecture method of classroom instruction:
a process by which the contents of the textbook of the instructor are 
transferred to the notebook of the student without passing through the head 
of either party.\rq\rq}
Here we undertake the endeavor of writing such an article
about crossed modules with an eye towards the past, this past
being ignored by the recent activity in this area.

In ancient cultures,
{\em symmetry\/} arose as  repetition of patterns.
The human being perceives such
symmetry as harmonious and beautiful proportion and balance
(music, art,  architecture, etc.).
The variety of patterns is  untold.
Symmetry enables us to structure this variety
by recognizing repetitions,
as Eurynome's dancing  structured chaos.
The symmetries of an object are encoded in transformations
that leave the object invariant.
In modern mathematics, abstracting from the formal properties
of such transformations led to the idea of a  group.
A typical example is  the group of 
symmetries of the solutions of
an equation, in modern mathematics
termed {\em  Galois group\/},
or each of the 17 plane symmetry groups.
The operations that form these groups 
were already known to the Greeks.
Besides  in mathematics, symmetry and groups play a major role
in physics, chemistry, engineering, materials science,
crystallography, meteorology, etc.
A group itself admits symmetries: the automorphisms of a group constitute
a group. 
The {\em crossed module\/} concept arises by
abstracting from the formal properties
which this pair enjoys.
For two groups $N$ and $Q$, 
the technology of crossed modules allows for a complete classification
of the groups $E$ having $N$ as a normal subgroup and  quotient group
$E/N$ isomorphic to $Q$. This settles the
{\em extension problem\/} for groups,
raised by Hoelder at the end of the 19'th century.
Schreier explored this problem
in terms of factor sets,
and Turing implicitly noticed that
it admits a solution in terms of crossed modules.
Given $N$ and $Q$ as symmetry groups,
interpreting an extension $E$ as a symmetry group is an interesting task,
as is,
given the extension group $E$ as a symmetry group,
interpreting  $N$ and $Q$ as symmetry groups.
Crossed modules occur
in mathematics under
various circumstances as a means to structure a collection
of mathematical patterns.
A  special example of a crossed module arises from
an ordinary module over a group.
While the concept of an ordinary module is
lingua franca in mathematics, this is not the case of crossed modules.

Beyond some basic algebra and algebraic topology, 
we assume the reader familiar with
some elementary category theory
(\lq\lq categorically thinking\rq\rq\  suffices). 
This article is addressed to the non-expert.
We keep sophisticated technology like group cohomology, homotopy theory, 
and algebraic number theory 
 at a minimum.
We write an injection as $\rightarrowtail$
and a surjection as $\twoheadrightarrow$.

The opener image of this article displays
the beginning of the first counterpoint
(contrapunctus)
 of the original printed version of 
J.~S.~Bach's \lq\lq Kunst der Fuge\rq\rq.
The order of the counterpoints (of the
second part thereof) had been lost
and, 100 years ago,  
Wolfgang Graeser, before enrolling as a mathematics student
at Berlin university, restored an  order 
(perhaps the original sequence) by means of 
symmetry considerations.
This order is the nowadays generally
accepted performance practice.

\section{Definition and basic examples}
A crossed module arises by abstraction from a structure
we are all familiar with when we run into a normal subgroup
of a group or into  the group of automorphisms of a group:
Denote the identity element of a group by $e$ and,
for a group $G$ and a $G$-(operator) group $H$ (a group $H$ together with 
an action
of $G$ on  $H$ from the left by automorphisms of $H$)
 we write the action as
$(x,y) \mapsto {}^xy$, for $x \in G$ and $y \in H$.
Consider two groups $C$ and $G$, 
an action of $G$ on $C$ from the left,
view the group $G$ as a $G$-group
with respect to conjugation,
and let 
$\partial \colon C \to G$ 
be a homomorphism of $G$-groups.
The triple $(C,G,\partial)$ constitutes
a {\em crossed module\/} if, furthermore,
for every pair $(x,y)$ of members of $C$,
\begin{equation}
 xyx^{-1}  ({}^{\partial x}y)^{-1}=e,
\label{peiff}
\end{equation}
that is, the members
$ xyx^{-1}$ and   ${}^{\partial x}y$ of $C$ coincide.
For a crossed module $(C,G,\partial)$,
the image $\partial(C)$
of $C$ in $G$ is a normal subgroup, 
the kernel $\ker(\partial)$ of $\partial$  is a central subgroup of $C$,
and the $G$-action on $C$ induces an action of the quotient group
$Q=G/(\partial( C))$ on $Z=\ker(\partial)$
turning $Z$ into a module over this group, and
it is common to refer to the resulting exact sequence
\begin{equation}
\xymatrix{
Z
\ar@{>->}[r]&
C
\ar[r]&
G
\ar@{->>}[r]&
Q
}
\label{crossed1}
\end{equation}
as a {\em crossed $2$-fold extension of 
$Q$ by $Z$\/}.
A morphism of crossed modules is defined in the obvious way.
Thus crossed modules constitute a category.
The terminology \lq crossed module\rq\  goes back to
\cite{MR30760}.
The identities \eqref{peiff} appear in
\cite{peiffone} and
 have come to be  known in the literature as
{\em Peiffer identities\/} \cite{MR47046}
(beware: not \lq Peiffer identity\rq\ as some of the
present day literature suggests).
There is a prehistory, however, 
\cite{zbMATH03012737, MR1557005}; 
in particular,
the Peiffer identities occur already in \cite{MR1557005}.
The reader 
may consult \cite[Section 3]{MR4322147} for details.

For a group $Q$ and a $Q$-module $M$,
the trivial homomorphism from $M$ to $Q$ is a crossed module structure
and, more generally, so is any $Q$-equivariant
homomorphism $\vartheta$  from $M$ to $Q$
such that the image $\vartheta(M)$ of $M$ in $Q$ acts trivially on $M$.
Relative to conjugation, 
the injection  of a normal subgroup
into the ambient group is manifestly a crossed module structure.
Also, it is immediate that the homomorphism
$\partial \colon G \to \Aut(G)$ 
from a group $G$ to its group $\Aut(G)$ of automorphisms
which sends a member of $G$ to the inner automorphism it defines
turns $(G,\Aut(G),\partial)$ into a crossed module.  It is common to
refer to the quotient group $\Out(G)=\Aut (G)/(\partial(G))$ 
as the group of {\em outer automorphisms\/}
of $G$.

\section{Identities among relations}
An \lq\lq identity among relations\rq\rq\ is for a presentation of a group
what is a \lq\lq syzygy among relations\rq\rq, as considered by Hilbert,
for a presentation of a module: 
Consider a presentation $\langle X;R\rangle$
of a group $Q$. That is to say, $X$ is a set of generators of $Q$
and $R$ a family of (reduced) words in $X$ and its inverses
such that the canonical epimorphism from
 the free group $F$ on $X$ to $Q$
has the normal closure $N_R$ of (the image of) $R$ in $F$ as its kernel.  
Then the $F$-conjugates of the 
images in $F$ of the
members
$r$ of $R$
generate $N_R$, 
that is, (with a slight abuse of the notation $R$,) the family
$R$  generates
 $N_R$ as an $F$-operator group:
With the notation ${}^y r = yry^{-1}$ ($r \in R$, $y \in F$), we
 can write any member $w$ of $N_R$ in the form
\begin{equation*}
w=\prod_{j=1}^m {}^{y_j}r_j^{\varepsilon_j},\ r_j \in R,\ y_j \in F,\ 
\varepsilon_j = \pm 1 ,
\end{equation*}
but $w$ does not determine such an expression uniquely. 
Thus the issue of understanding 
the structure of $N_R$ as an
$F$-operator group arises.
Heuristically, an {\em identity among relations\/}
is such a specified product that recovers the identity element $e$ of $F$.
For example, consider the presentation
\begin{equation}
\langle x,y; r,s,t\rangle,\ r = x^3,\, s = y^2,\, t = xyxy,
\label{examp2}
\end{equation}
of the symmetric group $S_3$ on three letters.
Straightforward verification shows that
\begin{gather}
t
s^{-1}
({}^{x^{-1}}t) 
({}^{x^{-1}}s^{-1})
({}^{x^{-1}y}r^{-1}) 
({}^{x^{-2}}t)
({}^{x^{-2}}s^{-1})
r^{-1}
\label{id1}
\end{gather}
is an  identity among the relations in \eqref{examp2}.
A standard procedure enables us to read off this identity
and others
from the following  prism shaped tesselated $2$-sphere:
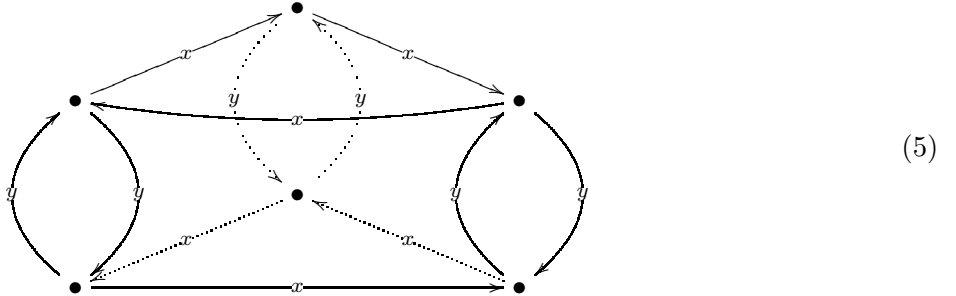
\begin{equation}
\begin{gathered}
\xymatrixcolsep{6pc}
\xymatrix{
  & \bullet 
\ar@{.>}@/_2pc/[dd]|-{y} 
 & & 
\\
 \bullet \ar@{->}
[ur]|-{x} 
\ar@{<-}@/_0.6pc/[rr] |-{x} 
 \ar@/^2pc/[dd]|-{y}
 & & \bullet 
\ar@/^2pc/[dd]|-{y} 
\ar@{<-}[ul]|-{x} 
\\
 & \bullet \ar@{.>}@/_2pc/[uu]|-{y}  
 & & 
 \\            
\bullet
\ar@/^2pc/[uu]|-{y}  
\ar@{->}[rr]|-{x}
\ar@{<.}[ur]|-{x}
 & & \bullet 
\ar@/^2pc/[uu]|-{y}  
\ar@{.>}[ul]|-{x}
}
\end{gathered}
\label{triangle}
\end{equation}
The reader will notice that reading 
 along the boundaries of the faces recovers the relators in 
\eqref{examp2}.
The projection, to a plane,
of this prism shaped $2$-sphere
with one of the triangles removed is
 \cite[Fig. 2 p.~155]{MR662431}. 
In Section 10 of  \cite{MR662431},  the reader 
can find precise methods to obtain 
such  identities from \lq\lq pictures\rq\rq, see in particular 
\cite[Fig. 12 p.~194]{MR662431} for the case at hand,
and in  \cite{MR675732} from diagrams etc.
The corresponding term in \cite{peiffone}
is \lq\lq Randwegaggregat\rq\rq.

To develop a formal understanding 
of the situation, let
$\widehat C_R$ be the free $F$-operator group on $R$.
The
kernel of the canonical epimorphism from the
free group on the (disjoint) union $X\cup R$ to $F$
 realizes $\widehat C_R$.
Let $\widehat \partial_R \colon \widehat C_R \to F$
denote the canonical homomorphism.
The members of the kernel of 
$\widehat \partial_R \colon \widehat C_R \to F$
are the {\em identities among
the relations\/} 
(or {\em among
the relators\/}) for the presentation  $\langle X;R\rangle$
 \cite{peiffone,reideone,MR662431};
Turing refers to them as \lq\lq relations between the relations\rq\rq\ 
 \cite[\S 2]{MR1557005}.

The {\em Peiffer elements\/}
$xyx^{-1}  ({}^{\widehat \partial_R x}y)^{-1}$, as 
$x$ and $y$ range over $\widehat C_R$, are
 identities that are
always present, independently of any particular presentation
under discussion.
Following \cite{reideone}, let
$C_R$ be the quotient group $\widehat C_R/P$ of
 $\widehat C_R$ 
modulo the  subgroup $P$ in $\widehat C_R$, necessarily normal,
which the Peiffer elements generate.
The Peiffer subgroup $P$ is an $F$-subgroup, whence 
the $F$-action on $\widehat C_R$ 
passes to an  $F$-action on $C_R$,
and the canonical homorphism
$\widehat \partial_R$ induces a homomorphism
$\partial_R \colon C_R \to F$
that turns $(C_R, F, \partial_R)$
into a crossed module.
In particular, the kernel $\pi= \ker(\partial_R)$,
being central in $C_R$, is an abelian group,  the $F$-action
on $\pi$ factors through a $Q$-module structure, and $\pi$
 parametrizes equivalence 
classes of \lq\lq essential\rq\rq\ 
or {\em non-trivial\/}
identities
associated to the presentation
$\langle X;R\rangle$ of  $Q$ \cite{reideone}.
We shall see below that there are interesting cases where
$\pi$ is trivial, that is, the Peiffer identities generate all identities.

Given a group $G$, a $G$-{\em crossed module\/}
is a $G$-group $C$ together with a homomorphism $\partial\colon C \to G$
of $G$-groups such that $(C,G,\partial)$ 
is a crossed module.
Given a set $R$ and a set map $\kappa \colon R \to G$
into a group $G$, the {\em free crossed\/}
$G$-{\em module on\/} $\kappa$
is the crossed module $(C_\kappa,G,\partial_\kappa)$
 enjoying the following property:
Given a $G$-crossed module $(C,G,\partial)$ and a set map
$\beta\colon R \to C$, there is a unique homomorphism
$\beta_\kappa \colon C_\kappa \to C$ of $G$-groups 
such that
$(\beta_\kappa,\Id) 
\colon 
(C_\kappa, G, \partial_\kappa) \to (C,G,\partial)$ 
is a morphism 
of crossed modules.
A standard construction shows that this free $G$-crossed
module always exists. 
By the universal property,  such a free $G$-crossed
module  is unique up to isomorphism,
whence it is appropriate to  use the definite article here.
The $F$-crossed module $(C_R, \partial_R)$ just constructed
from a presentation $\langle X; R\rangle$ of a group $Q$
is plainly the free $F$-crossed module on the injection
$R \to F$.
The resulting extension
$\xymatrix{
\pi
\ar@{>->}[r]&
C_R
\ar@{->>}[r]&
N_R}
$
of $F$-operator groups then displays the $F$-crossed module
$N_R$ 
as the quotient of a free
$F$-crossed module and thereby yields structural insight.
It is also common to
 refer to $(C_R, F, \partial_R)$
as the {\em free crossed module on\/}
$\langle X; R\rangle$.

\section{Group extensions and abstract kernels}

Given two groups $N$ and $Q$, the issue is to parametrize
the family of groups $G$ that contain $N$ as a normal subgroup
and have quotient $G/N$ isomorphic to $Q$, or equivalently,
in categorical terms,
the family of groups $G$ that fit into an exact sequence
of groups of the kind
\begin{equation}
\xymatrix{
& N \ar@{>->}[r]& G \ar@{->>}[r]& Q .
}
\label{ext1}
\end{equation}
Given an extension of the kind \eqref{ext1},
conjugation in $G$ induces a homomorphism 
$\varphi$ from $ Q$ to $\Out(N)$.
It is common to refer to a triple $(N,Q,\varphi)$
that consists of two groups $N$ and $Q$ together with
a homomorphism $\varphi \colon Q \to \Out (N)$
as an {\em abstract kernel\/}
or to the pair $(N,\varphi)$
as an {\em abstract $Q$-kernel\/}.
In \cite{zbMATH03012737} the terminology is \lq\lq Kollektivcharacter\rq\rq.
We have just seen that a group extension determines an abstract kernel.
Given two groups $N$ and $Q$ together with an 
abstract kernel structure $\varphi \colon Q \to \Out(N)$,
the extension problem consists in
{\em realizing\/} the abstract kernel, that is, in
parametrizing the extensions of the kind \eqref{ext1}
having $\varphi$ as its abstract kernel
provided  such an extension exists,
and the abstract kernel is then said to be {\em extendible\/}.

When $N$ is abelian, the group of outer automorphisms of $N$
amounts to the group $\Aut(N)$, an abstract $Q$-kernel structure 
on $N$ is equivalent to a $Q$-module structure on $N$, 
and the semi-direct product group
$N \rtimes Q$ shows that this abstract kernel is extendible.
When $N$ is non-abelian,
not every abstract $Q$-kernel $(N,\varphi)$ is extendible, however,
\cite[p.~415]{zbMATH03012737}. 
We shall shortly see that  non-extendible abstract kernels
abound in mathematical nature.

For two homomorphisms $G_1 \to Q$ and $G_2 \to Q$,
let  $G_1\times_Q G_2$ denote the {\em pullback\/} group,
that is, the subgroup of $G_1 \times G_2$ that consists of the pairs
$(x_1,x_2)$ that have the same image in $Q$.
 
For a crossed module $(C,G, \partial)$,
the action $G \to \Aut(C)$ of $G$ on $C$
plainly defines an abstract kernel
structure $\varphi\colon  G /(\partial(C)) \to \Out(C)$.
On the other hand, given the abstract kernel
$(N,Q,\varphi)$, the pullback group $G^\varphi = \Aut(N) \times_{\Out(N)} Q$ 
and  the canonical homomorphism
$\partial^\varphi\colon N \to G^\varphi$
which the crossed module structure 
$\partial\colon N \to \Aut(N)$
induces combine to the crossed module 
$(N,G^\varphi,\partial^\varphi)$,
and $\ker(\partial^\varphi)$ coincides with the center of $N$. 
In fact, this correspondence is a bijection
between abstract kernels and crossed modules
$(C,G,\partial)$
having $\ker(\partial)$ as the center of $C$.
The group $G^\varphi$ occurs in
\cite{zbMATH03012737}  as the \lq\lq Aufloesung des 
Kollektivcharakters\rq\rq\ $\varphi$.

We will now explain the solution of the extension problem:
Consider an abstract kernel
$(N,Q,\varphi)$,
and let $(N,G^\varphi,\partial^\varphi)$
be its associated
crossed module.
Let $\langle X; R\rangle$ be a presentation of the group $Q$
and, exploiting the freeness of $F$,  choose the set maps 
$\alpha\colon X \to G^\varphi$ and $\beta\colon R \to N$
in such a way that 
(i) the composite $F\to Q$ of the induced homomorphism 
$\alpha_F\colon F \to G^\varphi$
with the epimorphism from $G^\varphi$ to $Q$ coincides with the epimorphism
from $F$ to $Q$ and (ii)
the composite $R \to G^\varphi$  of $\alpha_F$ with the injection from $R$ to $F$
coincides with the composite of $\partial$ with $\beta$.
In view of the universal property of the free $F$-crossed module
$(C_R,\partial_R)$,
the set map $\beta$  induces a morphism
 $(\beta_R,\alpha_F)\colon (C_R,F,\partial_R) \to (N, G^\varphi,\partial^\varphi)$
of crossed modules
defined on the free crossed module $(C_R,F,\partial_R)$
on $\langle X; R\rangle$.
Let $\pi = \ker(\partial_R)$ and
display the resulting morphism of crossed modules as
\begin{equation}
\begin{gathered}
\xymatrix{
\pi
\ar[d]
\ar@{>->}[r]&
C_R
\ar[d]_{\beta_R}
\ar[r]^{\partial_R}&
F
\ar[d]^{\alpha_F}
\ar@{->>}[r]&
Q
\ar@{=}[d]
\phantom{.}
\\
Z
\ar@{>->}[r]&
N
\ar[r]\ar[r]_{\partial_\varphi}&
G^\varphi
\ar@{->>}[r]&
Q.
}
\end{gathered}
\label{diag4}
\end{equation}
(This is Diagram 8 in \cite[Section 3]{MR4322147}).
The following is a version of
\cite[Theorem 4 p.~356]{MR1557005} 
in modern language and terminology;
see \cite[Section 3]{MR4322147} for
 details.
Crossed modules, in particular free ones, exact sequences,
commutative diagrams,  pullbacks, etc.
were not at
Turing's disposal, however,
and he expressed his ideas in other ways.
\footnote{A. Turing and J.H.C. Whitehead
worked as WWII codebreakers at Bletchley Park (UK) but we shall never 
know whether they then discussed the idea of a crossed module.}

\smallskip
\noindent
{\bf {Theorem.}}
{\sl The abstract kernel
$(N,Q,\varphi)$ is
extendible if and only if, once
$\alpha$ has been chosen, the set map $\beta$ can be chosen in such a way that,
in Diagram {\rm \ref{diag4}},
the restriction of $\beta_R$ to $\pi$ is zero.
}

\begin{proof}
It is immediate that the condition is necessary.
To establish the converse,
suppose, $\alpha$ having been chosen, there is a choice of
 $\beta$ with the asserted property.
Let $N_R= \partial_R(C_R) \subseteq F$ 
denote the normal closure of $R$ in $F$.
From \eqref{diag4}, we deduce the morphism
\begin{equation}
\begin{gathered}
\xymatrix{
N_R
\ar[d]_{\widetilde \beta}
\ar@{>->}[r]^{\iota}&
F
\ar[d]^{\alpha_F}
\ar[r]&
Q
\ar@{=}[d]
\\
N
\ar[r]\ar[r]_{\partial_\varphi}&
G^\varphi
\ar[r]&
Q
}
\end{gathered}
\label{diag3}
\end{equation}
of crossed modules.
Since \eqref{diag3}
is a morphism of crossed modules,
in the semi-direct product group
$N \rtimes F$, the members of the kind
$(\widetilde \beta(y), \iota (y^{-1}))$ as $y$ ranges over $N_R$
constitute a  subgroup, necessarily normal.
The quotient group $E$ of 
$N \rtimes F$ by this normal subgroup
yields the  group extension
 $E$ of $Q$ by $N$  we seek.
\end{proof}

This  modern proof
of the Theorem 
is in \cite[Section 10]{MR576608}
(phrased in terms of the coequalizer of 
$\widetilde \beta$ and $\iota$).
The diligent reader will notice that the 
 proof is complete, i.e., no detail is left to the reader.
It is an instance of a common observation to the effect that
{\em mathematics consists in 
continuously improving
 notation and terminology\/}. 
Even nowadays there are textbooks that struggle to
explain the extension problem
and its solution in terms of lengthy and
 unilluminating
cocycle calculations, and it is hard to find
the above theorem,
 dug out as a result of Turing's only recently 
\cite{MR4322147}.\footnote{Turing would have slipped through an 
evaluation system
based on bibliographic metrics.}
Eilenberg-MacLane (quoting \cite{MR1557005}  but
apparently not understanding Turing's reasoning)
developed the obstruction for an abstract kernel to be extendible
in terms of the vanishing of the class of a
group cohomology 3-cocycle
of $Q$ with values in $Z$ \cite[Theorem 8.1]{MR0020996}.
The unlabeled arrow $\pi \to Z$ in Diagram \ref{diag4}
recovers this 3-cocycle, and the condition in the  Theorem 
says that the class of this $3$-cocycle
in the group cohomology group $\Ho^3(Q,Z)$
is zero.
Indeed, once a choice of $\alpha$ and $\beta$ has been made,
the vanishing of that $3$-cocycle is equivalent
to the unlabeled $Q$-module morphism $\pi \to Z$
in \eqref{diag4}
admitting an extension 
to a homomorphism
$\gamma \colon C_R \to Z$ of $F$-groups, and
setting
$\widetilde \beta_R(y) = \beta_R(y) \gamma(y)^{-1}$
and substituting 
$\widetilde \beta_R$ for $\beta_R$, we obtain
a diagram of the kind \eqref{diag4}
having $\pi \to Z$ zero.

To put flesh on the bones of the last remark, recall the 
{\em augmentation\/} map
$\varepsilon\colon \mathbb Z \Gamma \to \mathbb Z$ of a group 
$\Gamma$ 
defined on the integral group ring $\mathbb Z \Gamma$ of $\Gamma$ by
$\varepsilon (x) =1$ for $x \in \Gamma$ and
write the {\em augmentation ideal\/} 
$\ker(\varepsilon)$ as $\mathrm I \Gamma$.
The $F$-action on $C_R$
induces a $Q$-action on the
abelianized group $C_{R,\mathrm{ab}} = C_R/[C_R, C_R]$
in such a way that the set $R$ constitutes a set of free
generators; thus  $C_{R,\mathrm{ab}}$
is canically isomorphic to the free
$Q$-module $\mathbb Z Q[R]$ freely 
generated by $R$ (with a slight abuse of the notation $R$). 
The induced morphism
$\pi \to C_{R,\mathrm{ab}}$ of $Q$-modules is still injective. 
For $r \in R$ and $x \in X$, using the fact that 
 $\{x-1 \in \mathrm I F; x \in X\}$ is a family of free $F$-generators
of $\mathrm I F$,
define $r _x \in \mathbb Z F$ 
by the identity
$r-1= \sum  r _x (x-1)\in \mathrm I F$, and let $[r_x]\in \mathbb Z Q[X]$
denote the image.
Defining
 $d_2$ by $d_2(r) = \sum [r _x] x\in \mathbb Z Q[X]$ where 
$r$ ranges over $R$ and
 $x$ over $X$, and $d_1$ by $d_1(x)= [x]-1 \in \mathbb Z Q$ where
$[x] \in\mathbb Z Q$ denotes the image of $x$, as
 $x$ ranges  over $X$, renders 
the sequence
\begin{equation}
\xymatrixcolsep{1.2pc}
\xymatrix{
\pi
\ar@{>.>}[r]&
\mathbb Z Q[R]
\ar[r]^{d_2}&
\mathbb Z Q[X]
\ar[r]^{\phantom{aa}d_1}&
\mathbb Z Q
\ar@{.>>}[r]^{\varepsilon}& \mathbb Z}
\end{equation}
exact, and the 
three middle terms thereof 
together with the corresponding arrows 
constitute  
the beginning of a free resolution of $\mathbb Z$ in the category of
$Q$-modules. When we take $\langle X;R \rangle$ to be the 
{\em standard\/} presentation having $X$ the set of all members 
of $Q$ distinct from $e$, we obtain the first three terms of the
{\em standard resolution\/}, and
the composite of the 
$Q$-module epimorphism
$B_3 \to \pi$ defined on the fourth term $B_3$ of the standard resolution
with the above unlabeled arrow $\pi \to Z$ verbatim
recovers the $3$-cocycle in \cite[Theorem 8.1]{MR0020996}.

A {\em congruence\/} between two group extensions
of a group $Q$ by a group $N$ is a commutative diagram
\begin{equation}
\begin{gathered}
\xymatrix{
N
\ar@{=}[d]
\ar@{>->}[r]&
E_1
\ar[d]
\ar@{->>}[r]&
Q
\ar@{=}[d]
\\
N
\ar@{>->}[r]\ar[r]&
E_2
\ar@{->>}[r]&
Q
}
\end{gathered}
\label{diag12}
\end{equation}
with exact rows in the category of groups.
For a group $N$ with center $Z$, the multiplication map
restricts to a homomorphism $\mu\colon N \times Z \to N$.

\smallskip
\noindent
{\bf {Complement.}}
{\sl 
Let $(N,Q,\varphi)$ be an extendible abstract kernel and realize it
by the group extension
$N \stackrel{j_1}\rightarrowtail E_1 \twoheadrightarrow Q$.
The assignment to an extension $Z \stackrel {j}\rightarrowtail E \twoheadrightarrow Q$
of $Q$ by the center $Z$ of $N$ realizing the induced $Q$-module structure
on $Z$ of the group extension
$N \stackrel{j_2}\rightarrowtail E_2 \twoheadrightarrow Q$
which the commutative diagram
\begin{equation}
\begin{gathered}
\xymatrix{
N\times Z
\ar[d]_\mu
\ar@{>->}[r]^{(j_1,j)}&
E_1\times_Q E
\ar[d]
\ar@{->>}[r]&
Q
\ar@{=}[d]
\\
N
\ar@{>->}[r]\ar[r]&
E_2
\ar@{->>}[r]^{j_2}&
Q
}
\end{gathered}
\label{diag13}
\end{equation}
with exact rows
in the category of groups
characterizes
induces a faithful and transitive action of
$\Ho^2(Q,Z)$ on the congruence classes of extensions of $Q$ by $N$
realizing the abstract kernel $(N,Q,\varphi)$.
} 

As for the proof we note that,  as before, since 
the left-hand and middle vertical arrows in
\eqref{diag13} constitute a morphism of crossed modules,
in the group
$N \rtimes (E_1 \times _Q E)$,
the triples $(y_1y_2, j_1(y_1)^{-1}, j(y_2)^{-1})$, as $y_1$ ranges over 
$N$ and $y_2$ over $Z$,
form a  normal subgroup, and
the group $E_2$ arises as the quotient of 
$N \rtimes (E_1 \times _Q E)$ by this normal subgroup.
When $N$ is abelian, it coincides with its center,
Diagram \ref{diag13} displays the operation 
of \lq\lq Baer sum\rq\rq\ of two extensions with abelian kernel, 
and the complement recovers the fact that the
congruence classes of abelian extensions of $Q$ by $N$
constitute
an abelian group, the group structure being induced
by the operation of Baer sum
\cite{zbMATH03012737}, indeed,  
 the group cohomology group $\Ho^2(Q,N)$.

\section{Combinatorial group theory and low dimensional topology}
\label{comb}
Consider a space $Y$ and a subspace $X$ 
(satisfying suitable local properties,
$Y$ being a CW complex and $X$ a subcomplex would suffice), 
 and let $o$ be a base point in $X$.
The standard  action of the
fundamental group $\pi_1(X,o)$ (based homotopy classes of continuous maps from a
 circle to $X$ with a suitably 
defined composition law) on
the second relative homotopy group $\pi_2(Y,X,o)$
(based homotopy classes of continuous maps from a disk to $Y$ 
such that the boundary circle maps to $X$ with a suitably 
defined composition law) 
and the boundary map 
$\partial\colon  \pi_2(Y,X,o) \to\pi_1(X,o) $
turn $(\pi_2(Y,X,o),\partial)$ into a 
crossed $\pi_1(X,o)$-module
\cite{whitethr}.
See \cite{MR2841564} for a leisurely introduction to homotopy groups
building on an algebra of composition of cubes.
Suffice it to mention  that the higher homotopy groups 
of a space acquire an action of the fundamental group.
In \cite{whitethr},  J.H.C. Whitehead in particular proved that,
when the  space $Y$ arises
from $X$ by attaching $2$-cells,
the crossed $\pi_1(X,o)$-module
$(\pi_2(Y,X,o),\partial)$ is free on
the homotopy classes 
in $\pi_1(X,o)$
of the attaching maps of
 the $2$-cells.

The {\em Cayley graph\/}
of a group $Q$ with respect to a family $X$ 
of generators is the directed graph
having $Q$ as its set of vertices and,
for each pair $(x,y) \in X \times Q$,
an oriented edge joining $y$ to $xy$.
The oriented graph that underlies  \eqref{triangle}
is the Cayley graph of the group  $S_3$ with respect to the generators
$x$ and $y$ in \eqref{examp2}.
The {\em geometric realization\/}
$K=K\langle X;R\rangle$ 
of a presentation $\langle X;R\rangle$
of a group $Q$ 
is a $2$-dimensional CW complex
with a single zero cell $o$,
having $1$-cells in bijection with $X$ and $2$-cells
in bijection with $R$ in such a way that
the fundamental group 
$\pi_1(K^1,o)$
of the $1$-skeleton amounts to the free group
on $X$ and that the attaching maps of the $2$-cells
define, via the boundary map 
$\partial\colon \pi_2(K, K^1,o)\to \pi_1(K^1,o)$,
the members of $R$.
By construction, the fundamental group $\pi_1(K,o)$ of $K$ 
 is  isomorphic to
$Q$, the fundamental group $\pi_1(K^1,o)$ 
is canonically isomorphic to
the free group $F$ on the generators $X$
 and, by Whitehead's theorem,
the $\pi_1(K^1,o)$-crossed module  $(\pi_2(K, K^1,o),\partial)$
is free on the attaching maps of the $2$-cells
and hence canonically isomorphic to the
 $\pi_1(K^1,o)$-crossed module written above as $(C_R,\partial_R)$.
The $1$-skeleton $\overline K^1$ of the universal covering
space $\widetilde K$ of
$K$ then has fundamental group 
isomorphic to the normal closure $N_R$
of the relators in $F$ and, together with the appropriate orientation
of its edges,  then amounts to the Cayley graph of $Q$ with respect to $X$.
Every $2$-dimensional CW complex with a single $0$-cell is of 
this kind. Since the higher homotopy groups of $K^1$ are zero,
the long exact homotopy sequence of the pair 
$(K,K^1)$ reduces to the crossed $2$-fold extension
\begin{equation}
\xymatrixcolsep{1.1pc}
\xymatrix{
\pi_2(K,o)
\ar@{>->}[r]&
\pi_2(K, K^1,o)
\ar[r]^{\phantom{aa}\partial}&
\pi_1(K,o)
\ar@{->>}[r]&
Q.
}
\label{crossed2}
\end{equation}
Thus the group of \lq\lq essential identities\rq\rq\ 
among the relations $R$ 
amounts to the second homotopy group 
$\pi_2(K,o)$ of $K$ (based homotopy classes of continuous maps from a $2$-sphere to $K$,
with a suitably define group structure).
For illustration,
consider the $1$-skeleton of the
prism shaped tesselated $2$-sphere \eqref{triangle}.
For each $2$-gon and each $4$-gon, 
 attach two copies of the disk bounding
it and, likewise, for  each triangle,
attach three copies of the disk bounding it.
Thus  we obtain 
a prism shaped $2$-complex having, beyond the six vertices and twelve edges,
six faces bounded by $2$-gons,  six faces bounded by 
$4$-gons, and six faces bounded by triangles,
and this $2$-complex realizes the universal covering space $\widetilde K$
of the geometric realization $K$ of \eqref{examp2}.
By its very construction, it has eleven \lq\lq chambers\rq\rq\ 
(combinatorial $2$-spheres).
Since $\widetilde K$ is simply connected,
we conclude its second homotopy group is free abelian of rank eleven.
For the inexperienced reader we note that
choosing a maximal tree in $\overline K^1=\widetilde K^1$ 
and deforming this tree to point
is a homotopy equivalence, and the result is a bunch of eleven $2$-spheres,
necessarily having second homology group
free abelian of rank eleven.
Since $\widetilde K$ is simply connected, 
the Hurewicz map from $\pi_2(\widetilde K,o)$
to this homology group is an isomorphism, and so is
the homomorphism $\pi_2(\widetilde K,o) \to \pi_2(K,o)$
which the covering projection induces.
Hence
 $\pi_2(K,o)$ is free abelian of rank eleven.
Under the isomorphism 
$\pi_2(K,K^1,o) \to C_R$ of crossed $F$-modules,
the based homotopy class of the innermost chamber goes to the class
of the
 identity \eqref{id1} in $C_R$,
necessarily non-trivial, as are the identities that correspond 
to the other chambers. Another interesting piece of information
we extract form \eqref{triangle} is that  
$\pi_1(\overline K^1,o) \cong N_R \subseteq F$,
being the fundamental group of the Cayley graph $\overline K^1$ of $S_3$
is, as a group, freely
generated by any seven among the eight constituents of \eqref{id1}.

The paper \cite{peiffone}
arose out of a  combinatorial study of $3$-manifolds:
The present discussion applies to the $2$-skeleton $M^2$
of a  cell decomposition of a $3$-manifold $M$,
and attaching $3$-cells to build the $3$-manifold under discussion
\lq\lq kills\rq\rq\  some of the essential identities associated with
$M^2$:
The commutative diagram
\begin{equation}
\begin{gathered}
\xymatrixcolsep{1.1pc}
\xymatrix{
\pi_2(M^2,o)
\ar@{->>}[d]
\ar@{>->}[r]&
\pi_2(M^2,M^1,o)
\ar@{->>}[d]
\ar[r]^{\phantom{aaa}\partial}&
\pi_1(M^1,o)
\ar@{=}[d]
\ar@{->>}[r]&
\pi_1(M,o)
\ar@{=}[d]
\phantom{.}
\\
\pi_2(M,o)
\ar@{>->}[r]&
\pi_2(M,M^1,o)
\ar[r]\ar[r]_{\phantom{aaaa}\partial_M}&
\pi_1(M^1,o)
\ar@{->>}[r]&
\pi_1(M,o)
}
\end{gathered}
\label{diag5}
\end{equation}
with exact rows
displays how the $\pi_1(M^1,o)$-crossed module
$(\pi_2(M,M^1,o),\partial_M)$ arises from
the free $\pi_1(M^1,o)$-crossed module 
$(\pi_2(M^2,M^1,o),\partial)$.
For these and related issues, see, e.g., 
\cite{MR662431} and the literature there.
For illustration,
consider a finite subgroup $Q$ of $\mathrm {SU}(2)$ such as, e.g.,
$Q=C_n$, the cyclic group of order $n\geq 2$,
or $Q$ the quaternion group of order eight.
The orbit space $M=\mathrm {SU}(2)/Q$, a {\em lens space\/} 
when $Q=C_n$,
is a $3$-manifold having $\pi_2(M^2,o)$ non-trivial but  $\pi_2(M,o)$
trivial since the $3$-sphere has trivial second homotopy group.
Via Papakyriakopoulos' {\em sphere theorem\/},
the  second homotopy group of a general 
$3$-manifold  being non-trivial is equivalent 
to the manifold being  geometrically  splittable, 
similarly to what we hint at for
tame links below. Thus a  $3$-manifold is geometrically splittable if 
and only if attaching the $3$-cells of a cell decomposition does not kill
all the essential identities arising from its $2$-skeleton.

A $2$-complex $K$ is said to be {\em aspherical\/} when its 
degree $\geq 2$  homotopy groups  are trivial.
This is equivalent to the second homotopy group  being trivial.
Lyndon's identity theorem \cite{MR47046}
implies that, when 
$K$ has a single $2$-cell, 
the $2$-complex $K$ is aspherical if and only if
the relator
$r$ arising as the boundary of the $2$-cell
is not a proper power in $F = \pi_1(K,o)$.
Thus the  $F$-crossed module structure
of 
the normal closure $N_{\{r\}}$ of $\{r\}$ in $F$ is free in this case.
For example, a closed surface distinct from the $2$-sphere
is aspherical, but this is an immediate consequence of
the universal covering being the $2$-plane.
In the same vein, the exterior of a tame link  in the $3$-sphere
is homotopy equivalent to a  
$2$-complex--the operation of \lq\lq squeezing\rq\rq\ 
the $3$-cells of a cell decomposition 
achieves this--, and the
link is geometrically unsplittable if and only 
if the exterior and hence the $2$-complex is aspherical. 
Here \lq\lq tame\rq\rq\ means that the link belongs to a cellular
subdivision.
Consequently,
for the Wirtinger presentation $\langle X;R\rangle$ 
\cite[Fig. 3 \S 9 p.~183]{MR662431} of the link group,
as a crossed $F$-module,
the normal closure $N_R$ of the relators $R$ in $F$ is free
if and only if, by the already  quoted result of Papakyriakopoulos,
 the link is geometrically unsplittable
\cite[Theorem (P) \S 9 p.~183]{MR662431}.
In particular, a tame knot is a geometrically unsplittable link.
In \cite{whitethr}, Whitehead raised the issue,
still unsettled,
 whether
any subcomplex of an aspherical $2$-complex is itself aspherical.
This is equivalent to asking whether, for a presentation
 $\langle X;R\rangle$ of a group,
when the normal closure $N_R$ of $R$ in $F$
is free as an $F$-crossed module,
this is still true of the normal closure $N_{\widetilde R}$ of 
a subset $\widetilde R$ 
of $R$. See \cite[Section 9 p.~181 ff.]{MR662431} for more details 
and literature on partial results.
In \cite{whitethr}, in the proof of the freeness of the crossed module
arising from attaching $2$-cells to a space,
Whitehead interpreted the Peiffer identities as Wirtinger relations
associated to the link arising from a null homotopy;
see \cite[Section 10 p.~187 ff.]{MR662431} for details and more references.

\section{Interpretation of the third group cohomology group}

The notion of congruence for group extensions extends to 
crossed $2$-fold extensions,
and congruence  classes of crossed $2$-fold extensions
of the kind \eqref{crossed1}, together with a suitably 
defined  operation of composition arising from a generalized Baer sum,
constitute the third group cohomology group $\Ho^3(Q,Z)$.
Under this interpretation, the crossed $2$-fold extension
associated to the crossed module $(Z,Q,0)$
represents the identity element.
See \cite{histnote} for the history of this interpretation;
it parallels the result in \cite{MR0020996} saying that
suitably defined equivalence classes of abstract $Q$-kernels
having $Z$ as its center 
are in bijection with the members of 
 $\Ho^3(Q,Z)$.
Under this interpretation, the zero class corresponds to the 
extendible $Q$-kernels.
This is essentially Turing's theorem in a new guise.
Thus  the crossed $2$-fold extension 
\eqref{crossed1} associated to
a $G$-crossed module $(C,\partial)$ defines a 
{\em characteristic class\/}
in $\Ho^3(Q,Z)$ with $Q= G/\partial(C)$.
In particular, such a
$G$-crossed module
having $\ker(\partial)$ equal to the center $Z$ of $C$
defines an abstract $Q$-kernel,
and this $Q$-kernel is extendible if and only if
its characteristic class is zero.
The characteristic class in $\Ho^3(\pi_1(K,o),\pi_2(K,o))$
of the $\pi_1(K^1,o)$-crossed module
$(\pi_2(K,K^1,o),\partial)$ associated with a CW complex $K$
recovers the (first) $k$- (or {\em Postnikov\/}) invariant;
when $K$ is $2$-dimensional, 
$\pi_1(K,o)$, $\pi_2(K,o)$ and this $k$-invariant
determine the homotopy type of $K$.
A CW complex with non-trivial fundamental group and non-trivial 
second homotopy group
 typically has non-zero $k$-invariant
but, to understand the present discussion,
 there is no need to know anything about
$k$-invariants beyond the fact that the crossed $2$-fold extension
associated with the corresponding crossed module represents it.
The
crossed module associated with the
geometric realization of \eqref{examp2}
provides an explicit example of a non-trivial  $k$-invariant, see below.
Here is an even more elementary example, with explicit verification
of the non-triviality of its $k$-invariant:
The
free crossed module $(C_{\{r\}},C_x,\partial)$
associated with the presentation $\langle x,r \rangle$ 
with $r = x^n$ of the finite cyclic group $C_n$ of order $n \geq 2$
has $C_x$ the free cyclic group generated by $x$
and $C_{\{r\}}$ the free abelian $C_n$-group
which $r$ generates, 
equivalently, the free $\mathbb Z C_n$-module which $r$ generates
when we use additive notation,
the action of $C_x$ on $C_{\{r\}}$ is the composite of the projection
$C_x \to C_n$ with the action coming from the $C_n$-group structure
on $C_{\{r\}}$, and
$\partial$ sends ${}^{x^k}r$, $0 \leq k \leq n-1$, to $x^n \in C_x$.
The $n$ identities
\begin{equation*}
i_1 ={}^x \negthinspace r r^{-1}, 
i_2 ={}^x \negthinspace i_1 ={}^{x^2} r ({}^x\negthinspace r^{-1}), \ldots ,
i_n ={}^{x^{n-1}}  i_1=r ({}^{x^{n-1}} r^{-1}) ,
\end{equation*}
generate  the kernel $\pi=\ker(\partial)$
as an abelian group, 
 subject to the relation
$i_1 i_2 \ldots i_n = 1$.
Thus $i_1$ generates $\pi$ as an abelian
 $C_n$-group, equivalently, as a $\mathbb C_n$-module when we write
$\pi$ additively.
The experienced reader will recognize $\pi$ as being,
as a $C_n$-module, isomorphic to
the augmentation ideal 
$\mathrm I C_n$ of $C_n$.
Sending each $i_j$, for $1 \leq j \leq n$, to the generator of 
(a copy of)
$C_n$
defines an epimorphism $\pi \to C_n$; indeed, this
is the epimorphism that arises by dividing out the
$C_n$-action on $\pi$.
Let $(\widehat C_{\{r\}},\widehat \partial)$ denote the
$C_x$-crossed module which requiring the diagram
\begin{equation}
\begin{gathered}
\xymatrix{
\pi
\ar@{->>}[d]
\ar@{>->}[r]&
C_{\{r\}}
\ar@{->>}[d]
\ar[r]^{\phantom{a}\partial}&
C_x
\ar@{=}[d]
\ar@{->>}[r]&
C_n
\ar@{=}[d]
\phantom{.}
\\
C_n
\ar@{>->}[r]&
\widehat C_{\{r\}}
\ar[r]\ar[r]_{\phantom{a}\widehat\partial}&
C_x
\ar@{->>}[r]&
C_n
}
\end{gathered}
\label{diag6}
\end{equation}
with exact rows to be commutative characterizes.
Ler $C_r \subseteq C_x$ denote the free cyclic subgroup which $r =x^n$
generates.
As an abelian group,
$\widehat C_{\{r\}} \cong C_n \times C_r$ and, with the notation
$u$ for the generator of the copy of $C_n$, the rules
${}^xr = ur$ and ${}^xu =u$
characterize the $C_x$-group structure.
Let $v$ denote a generator of the  cyclic group $C_{n^2}$ of
order $n^2$. 
With respect to the trivial actions,
the homomorphism
${\,\cdot\,}^n\colon C_{n^2} \to C_{n^2}$
which sends $v$ to $v^n$
defines a 
$ C_{n^2}$-crossed module structure on $ C_{n^2}$.
Sending
$u$ to $ v^n$, $r$ to $v$ and $x$ to $v$ we obtain
a congruence morphism
\begin{equation}
\begin{gathered}
\xymatrix{
C_n
\ar@{=}[d]
\ar@{>->}[r]&
\widehat C_{\{r\}}
\ar@{->>}[d]
\ar[r]^{\phantom{a}\widehat \partial}&
C_x
\ar[d]
\ar@{->>}[r]&
C_n
\ar@{=}[d]
\phantom{.}
\\
C_n
\ar@{>->}[r]&
 C_{n^2}
\ar[r]\ar[r]_{\phantom{a}{\,\cdot\,}^n}&
C_{n^2}
\ar@{->>}[r]&
C_n
}
\end{gathered}
\label{diag7}
\end{equation}
of crossed $2$-fold extensions.
The upper row of \eqref{diag6} represents the $k$-invariant 
in $\Ho^3(C_n, \pi)$ of
the geometric realization of the presentation
$\langle x;r\rangle$ of the group $C_n$ while the lower row of \eqref{diag7}
represents a generator of $\Ho^3(C_n,C_n) \cong C_n$.
Diagram \eqref{diag6} says that,  under the induced map
 $\Ho^3(C_n,\pi) \to\Ho^3(C_n,C_n)$,
that $k$-invariant goes to a generator of a cyclic group of order $n \geq 2$.
Hence that $k$-invariant is necessarily non-trivial.
It is also worthwhile noting that
the left-hand copy of $C_n$ 
in \eqref{diag7} amounts to
 the third
group homology group $\Ho_3(C_n)$ of $C_n$ and that
suitably exploiting
 the {\em Yoneda\/} 
interpretation of 
the traditional definition of
$ \Ho^3(C_n,C_n)$ as
$\mathrm{Ext}^3_{C_n}(\mathbb Z, C_n)$
identifies the class of the bottom rom
in \eqref{diag7} with the corresponding member of
 $\mathrm{Ext}^3_{C_n}(\mathbb Z, C_n)$.
Restricting the crossed $2$-fold extension  associated with the
geometric realization of \eqref{examp2}
to any of the cyclic subgroups of $S_3$ and
playing a bit with the data, one can also show that
the  crossed module associated with the
geometric realization of \eqref{examp2} has non-zero $k$-invariant.
Thus  crossed modules having non-zero characteristic class 
and in particular non-extendible abstract kernels
abound.

\section{Higher group cohomology groups}
Suitably extended, the interpretation in terms of crossed $2$-fold extensions
leads, 
for $n \geq 1$, to an interpretation of
the  group cohomology group $\Ho^{n+1}(Q,Z)$
in terms of \lq\lq crossed $n$-fold extensions\rq\rq.
See \cite{histnote} for the history of this interpretation.
For example, 
the crossed $n$-fold extension associated with a cell decomposition of
an $n$-dimensional CW-complex $X$
with non-trivial fundamental group and trivial homotopy groups
$\pi_j(X)$ for $2 \leq j < n$ represents 
the first non-zero $k$-invariant in $\Ho^{n+1}(X,\pi_n(X))$
of $X$ \cite{MR576609}. For illustration,
as in Section \ref{comb},
consider the orbit space $M=\mathrm {SU}(2)/Q$
for
a (non-trivial) finite subgroup $Q$ of $\mathrm {SU}(2)$.
We will now use, without further explanation, some classical material
which the reader can find in standard textbooks.
The $2$-skeleton of a suitable 
cell decomposition of $M$
yields the geometric realization of 
the presentation of $Q$ 
resulting from the cell decomposition,
 the geometric realization of 
the presentation $\langle x;r\rangle$ of $C_n$, with $r= x^n$,
when $Q=C_n$.
For
a general finite subgroup $Q$ of $\mathrm {SU}(2)$,
the exact homotopy sequence 
of the pair $(M,M^2)$, necessarily one of $Q$-modules, takes the form
\begin{equation}
\ldots  \pi_3(M^2,o)\to
\pi_3(M,o)
\to
\pi_3(M,M^2,o)
\to
\pi_2(M^2,o)\twoheadrightarrow
\pi_2(M,o) .
\label{form9}
\end{equation}
With respect to  the epimorphism
from $\pi_1(M^1,o)$ to $Q$,
\eqref{form9}
becomes a sequence of $\pi_1(M^1,o)$-modules.
The composite
\begin{equation}
\pi_3(M,M^2,o)
\to\pi_2(M^2,M^1,o) \to \pi_2(M^2,M^1,o)_{\mathrm{ab}}
\end{equation}
of
$\pi_3(M,M^2,o)
\to
\pi_2(M^2,o)
$
with the injection $\pi_2(M^2,o) \to \pi_2(M^2,M^1,o)$
and, thereafter, with abelianization,
amounts to the boundary operator $C_3(S^3) \to C_2(S^3)$
of the $Q$-equivariant 
cellular chain complex $C_*(S^3)$
 of the (cellularly decomposed) 
$3$-sphere $S^3$ that underlies $\mathrm {SU}(2)$, and this boundary operator
has the homology group $\Ho_3(S^3)$ as its kernel.
Using the Hurewicz isomorphism
$\pi_3(S^3,\widehat o) \to \Ho_3(S^3)$ and the covering projection isomorphism
$\pi_3(S^3,\widehat o) \to\pi_3(M,o)$ (with the notation
$\widehat o$ for a pre-image in $S^3$ of the base point $o$ of $M$), we deduce
that the $Q$-morphism 
$\pi_3(M,o)
\to
\pi_3(M,M^2,o)$ in \eqref{form9}
is injective. (Beware: We must be circumspect at this point
since $\pi_3(M^2,o)$ is non-trivial when $Q$ is not the  trival group, 
and we cannot
naively deduce that injectivity from the exactness of 
\eqref{form9}.)
Since the second  homotopy 
group of the $3$-sphere $S^3$ 
is trivial,
so is $\pi_2(M,o)$.
Consequently the $Q$-morphism 
$
\pi_3(M,M^2,o) \to \pi_2(M^2,o)$ in \eqref{form9}
is surjective.
Hence splicing \eqref{crossed2}, with $M^2$ substituted for $K$,
and
\eqref{form9}
yields
the crossed $3$-fold extension
\begin{equation}
\xymatrixcolsep{0.7pc}
\xymatrix{
\pi_3(M,o)
\ar@{>->}[r]&
\pi_3(M,M^2,o)
\ar[r]&
\pi_2(M^2,M^1,o)
\ar[r]^{\phantom{aaa}\partial}&
\pi_1(M^1,o)
\ar@{->>}[r]&
Q 
}
\label{crossed13}
\end{equation}
 of $Q$ by  the free cyclic group $\pi_3(M,o) \cong \Ho_3(S^3)$ 
(as $Q$-modules),
and the $Q$-action on 
 $\pi_3(M,o)$ is trivial since 
 this action amounts to the induced $Q$-action
on $\Ho_3(S^3)$, necessarily trivial.
The group $\Ho^4(Q,\pi_3(M,o))$
is cyclic of order $|Q|$ (the number of elements of $Q$), and
the crossed $3$-fold extension \eqref{crossed13}
represents a generator of 
$\Ho^4(Q,\pi_3(M,o))$
and thence the first non-zero $k$-invariant of $M$.
This $k$-invariant,
the fundamental group
$\pi_1(M,o) \cong Q$, and the third homotopy group $\pi_3(M,o)$ 
determine the homotopy type of $M$.
Furthermore, the group $Q$ has periodic cohomology, of period $2$
when $Q$ is cyclic and of period $4$ otherwise.
Thus, the operation of
cup product with the class of \eqref{crossed13}
induces isomorphisms $\Ho^s(Q,\,\cdot\,) \to \Ho^{s+4}(Q,\,\cdot\,)$,
for $s \geq 1$, for any integer $s$ when we interpret
the notation  $\Ho$ as {\em Tate\/} cohomology.
When $Q$ is a cyclic group $C_n$ of order $n>1$,
 the extension 
$\xymatrix{\pi_3(M,o)\ar@{>->}[r]&\pi_3(M,M^2,o)\ar@{->>}[r]&\pi = \ker(\partial)
}$
of $C_n$-modules
amounts to the familiar extension
$\xymatrix{\mathbb Z \ar@{>->}[r]&\mathbb Z C_n\ar@{->>}[r]& 
\mathrm I C_n
}$
of $C_n$-modules 
arising from the standard small 
free resolution of $\mathbb Z$
in the category of $C_n$-modules, and it is immediate that
 the  $C_n$-module structure
on $\pi_3(M,o)$ is trivial.

\section{Generalization}

The description of group cohomology in terms of crossed $n$-fold extensions
($n \geq 1$)
is susceptible to generalizations where
cocycles are not necessarily available.
For example, 
for an extension of a topological group $G$ by a continuous $G$-module
whose underlying bundle is non-trivial,
(global) continuous cocycles are not available.
See \cite[Section 3]{MR4322147}
and the literature there
for more situations where this happens.

\section{Lie algebra crossed modules}
The axiom \eqref{peiff} 
makes perfect sense for Lie algebras, and
the interpretation of the $(n+1)$th Lie algebra cohomology
group
in terms of crossed $n$-fold extensions ($n \geq 1$) is available.
Also the abstract kernel concept extends to Lie algebras
in an obvious manner, as does the equivalence
between abstract kernels and crossed modules
with the central kernel constraint explained above,
and there is an analogue of Turing's theorem.
Indeed, much of the above material carries over  to Lie algebras.
See \cite[Section 5]{MR4322147}
and the literature there for details.

\section{Normality of a non-commutative algebra over its center}
Crossed modules arise in Galois theory.
We will now  briefly delve into  this.
See \cite[Section 4]{MR4322147}
for the history 
and 
\cite{MR3769367}
for a more complete account and references:

Let $S$ be a commutative ring and
$A$ an  $S$-algebra having $S$ as its center.
Let 
$Q$ be a group of operators on $S$.
The development of a Galois theory for such algebras leads to the 
following question:
{\em Does every automorphism in $Q$
extend to an automorphism of\/} $A$?
The 
algebra $A$ is said to be $Q$-{\em normal\/}
when this happens to be the case.
We formalize the situation as follows:

Denote by $\Aut(A)$ the group of ring
automorphisms
of $A$ and by $\mathrm U(A)$ its group of units. 
The 
obvious homomorphism
$\partial \colon \mathrm U(A) \to \Aut (A)$
assigns to a unit of $A$ the associated inner automorphism of $A$, 
the obvious action of $\Aut(A)$
on $\mathrm U(A)$ turns the triple $(\mathrm U(A),\Aut(A), \partial )$ into a crossed module, and
$\ker (\partial) = \mathrm U(S)$, the group of units of $S$.
Write 
$\Out(A) = \Aut(A)/(\partial \mathrm U(A)) $.
Each inner automorphism of $A$ leaves $S$ elementwise fixed whence
 the  
restriction map $\Aut(A) \to \Aut(S)$ induces a homomorphism 
${\res\colon \Out(A) \to \Aut(S)}$.
Let 
$Q$ be a group and
$\kappaQ \colon Q\to \Aut(S)$ an action of  $Q$ on
$S$ by ring automorphisms.
Define a $Q$-{\em normal structure\/} on 
the central $S$-algebra $A$ 
{\em relative to the given action 
$\kappaQ \colon Q \to \Aut(S)$
of $Q$ on\/} $S$
to be a
homomorphism $\sigma\colon Q \to \Out(A)$ that lifts the action 
$\kappaQ \colon Q \to \Aut(S)$ 
of $Q$ on $S$
in the sense that the composite of
$\sigma$ with 
 $\res\colon \Out(A) \to \Aut(S)$ coincides with $\kappaQ$.
A $Q$-{\em normal $S$-algebra\/} is, then, a central $S$-algebra
$A$ together with a $Q$-normal structure $\sigma\colon Q \to \Out(A)$.

Let $(A, \sigma )$ be a $Q$-normal $S$-algebra,
let $G^\sigma ={\Aut (A) \times_{\Out (A)}Q}$, and let $G^\sigma $
act on $\mathrm U(A)$ 
via the canonical homomorphism from 
$G^\sigma $ to $\Aut(A)$.
The homomorphism
$\partial ^\sigma \colon \mathrm U(A) \to G^\sigma$ 
which  $\partial \colon \mathrm U(A) \to \Aut(A)$
induces turns
$(\mathrm U(A),G^\sigma, \partial^\sigma)$
into a crossed module, and
the crossed 2-fold extension
\begin{equation}
\mathrm e_{(A,\sigma)} \colon \mathrm U(S) 
\rightarrowtail \mathrm U(A) 
\stackrel{\partial^\sigma}\to G^\sigma 
\twoheadrightarrow Q
\label{pb11}
\end{equation}
represents a class
$[\mathrm e_{(A,\sigma)}] \in \mathrm H^3 (Q, \mathrm U(S))$,
the {\em Teichmueller class\/}
of $(A,\sigma)$.
For the special case where $S$ is a field,
a cocycle description of this class 
(independently of any crossed module)
is in \cite{MR0002858}.\footnote{Teichmueller
worked as WWII codebreaker for the high command of the German army.}

Define a $Q$-{\em equivariant\/} $S$-algebra 
to be a central $S$-algebra
$A$ together with a homomorphism $\rho\colon Q \to \Aut(A)$
that induces the $Q$-action $\kappa$ on $S$.
For example, 
let $R=S^Q$, the subring of $Q$-invariants in $S$;
the $S$-algebra 
$A=B \otimes_R S$ for some central $R$-algebra $B$
plainly admits a canonical
$Q$-equivariant structure.
Consider  a $Q$-normal $S$-algebra $(A,\sigma)$.
We can then ask, as did Teichmueller in the situation 
he considered,
 whether 
the $Q$-action on $S$ lifts to a $Q$-equivariant
structure.
When such a lift exists,
it induces a congruence between 
$\mathrm e_{(A,\sigma)}$ and the
crossed $2$-fold extension
arising from the crossed module $(Z,Q,0)$,
and hence
 the class  $[\mathrm e_{(A,\sigma)}]\in\mathrm H^3 (Q, \mathrm U(S))$ 
is zero. As for the converse, 
let
 $\mathrm M_I(A)$ denote the $(I \times I)$ matrix algebra
over $A$ for an index family $I$; when $I$ is not finite, 
we interpret $\mathrm M_I(A)$ 
as being
the endomorphism ring
of $\oplus_I A^{\mathrm{op}}$. The algebra $\mathrm M_I(A)$ 
is again a central $S$-algebra.
It is obvious that an automorphism of $A$ yields one of $\mathrm M_I(A)$
in a unique way, and the obvious map $A \to \mathrm M_I(A)$ is a ring homomorphism.
Hence a $Q$-normal structure $\sigma \colon Q \to \Out(A)$ on $A$
determines one on $\mathrm M_I(A)$, 
and we
denote this structure by $\sigma_I \colon Q \to \Out(\mathrm M_I(A))$.
By \cite[Theorem 6.1]{MR3769365},
the Teich\-m\"uller class
of a $Q$-normal $S$-algebra $(A,\sigma )$
is zero if and only if, for $I = Q$, the $Q$-normal structure $\sigma_I$
on the matrix algebra $\mathrm M_I(A)$ comes from an  equivariant one.
Thus the  class $[\mathrm e_{(A,\sigma)}]\in \Ho^3(Q, \mathrm U(S))$ is
the obstruction 
for the $Q$-normal algebra $(A,\sigma)$ to be equivalent to
a $Q$-equivariant one in the sense just explained.
In general, we cannot have 
$(A,\sigma)$ itself to be equivariant.
See  \eqref{fit} below for a special case.

Suppose $S$ is a field. Then we are running into ordinary Galois theory.
Here is a family of  explicit examples
of a $Q$-normal algebra having non-trivial Teichmueller class:
Consider a field $K$ and let
$L=K(\zeta)$  be  a normal  extension
having Galois group $N$ cyclic of order $n\geq 2$ (say).
Let $\tau$ denote a generator of $N$,  
let $\eta \in \mathrm U(K)$,
and 
consider the cyclic central simple $K$-algebra
$D(\tau,\eta)$ generated by 
$L=K(\zeta)$ and some (indeterminate)
$u$ subject to the relations
\begin{equation}
u \lambda = {}^\tau\negthinspace\lambda u,\ u^n = \eta,\ \lambda \in 
L=K(\zeta).
\label{rels}
\end{equation}
In $D(\tau,\eta)$, the member  $u$ is a unit having inverse
$u^{-1}= u^{n-1}\eta^{-1}$ and,
for $\lambda \in L$,
we get
 $u \lambda u^{-1}= {}^\tau\negthinspace\lambda$, that is,
the action of the Galois group $N$ on $L$
extends to the inner automorphism of  $D(\tau,\eta)$
which $u$ determines. 
The algebra 
$D(\tau,\eta)$ is
 a crossed product of $N$ with 
$L$ relative  to the $\mathrm U(L)$-valued  
2-cocycle of $N$
determined by $\eta$
but this fact need not concern us here.
The field $L$ is a maximal commutative subalgebra
of $D(\tau,\eta)$.
The capital $D$ serves as a mnemonic for the fact that
Dickson explored such algebras.

Distinct choices of $\eta \in \mathrm U(K)$ may lead to the \lq same\rq\ 
 algebra 
of the kind $D(\tau,\eta)$: 
The assignment to $\vartheta \in \mathrm U(L)$ of 
$\prod_{j=0}^{n-1} {}^{\tau^j}\negthinspace\vartheta$
defines the classical  {\em norm map\/} 
$\nu \colon \mathrm U(L)\to \mathrm U(K)$. 
Let $\vartheta \in \mathrm U(L)$.
Then
$u \mapsto u_\vartheta = \vartheta u$ induces an 
isomorphism $\alpha_\vartheta \colon 
D(\tau,\eta) \to D(\tau,\nu(\vartheta)\eta) $
which restricts to the identity of $L$, 
an automorphism
$\alpha_\vartheta$
of $D(\tau,\eta)$
if and only if $\nu(\vartheta)=1$.
Furthermore, for 
$\eta =\nu(\vartheta)$
with $\vartheta \in \mathrm U(L)$,
the algebra $D(\tau,\eta)$
comes down to  the algebra of $(n \times n)$-matrices over $K$.
Thus  the cokernel $\mathrm{coker}(\nu)$ of the norm map
$\nu \colon \mathrm U(L)\to \mathrm U(K)$ 
parametrizes  classes of algebras of the kind $D(\tau,\eta)$ 
such that
 $D(\tau,\eta_1)$ and  $D(\tau,\eta_2)$
belong to the same class if and only if
an isomorphism which restricts to the identity of $L$
carries $D(\tau,\eta_1)$ to $D(\tau,\eta_2)$.
The expert will recognize that
$\mathrm{coker}(\nu)$ amounts to $\Ho^2(N,\mathrm U(L))$
and $\mathrm{ker}(\nu)$ to the group of
multiplicatively written 
$\mathrm U(L)$-valued $1$-cocycles of $N$.

Let
$\Aut_L(D(\tau,\eta))$ denote the  group of automorphisms of 
$D(\tau,\eta)$ that restrict to an automorphism of $L|\fiel $.
For $\chi \in \mathrm U(L)$, the member
$\vartheta_\chi = {}^\tau\negthinspace \chi \chi^{-1}$
of $\mathrm U(L)$
lies in the kernel of $\nu$.
{\em The assignment to  $\chi$ of  
$\alpha_{\vartheta_\chi}\in \Aut_L(D(\tau,\eta))$
induces an embedding
of
$\mathrm U(L) /\mathrm U(K)$ into $\Aut_L(D(\tau,\eta))$
onto the subgroup of automorphisms
that restrict to the identity of $L$.}
Indeed,
every automorphism $\alpha$ of 
$D(\tau,\eta)$ that restricts to the identity
of $L$
necessarily satisfies the identity
\begin{align*}
\alpha(u) u^{-1} \lambda u \alpha(u)^{-1} &= 
\alpha(u)  ({}^{\tau^{-1}}\lambda) \alpha(u)^{-1}
= 
\alpha(u ({}^{\tau^{-1}}\lambda) u^{-1})= \alpha(\lambda) = \lambda
\end{align*}
for every $\lambda \in L$,
and this implies $\alpha(u) u^{-1} \in \mathrm U(L)$
since $L$ is a maximal commutative subalgebra
of $D(\tau,\eta)$; then $\vartheta =\alpha(u) u^{-1}$
belongs to the kernel of $\nu$.
By Hilbert's \lq\lq Satz 90\rq\rq,   
every member $\vartheta$ of $\ker(\nu)$ 
is of the kind $ \vartheta_\chi$, 
for some $\chi \in \mathrm U(L)$.

Let $Q$ be a finite group of operators on $K$ and let
$\fiel = K^Q$, so that $K|\fiel$ is a Galois extension.
Suppose that, furthermore,
$L|\fiel$ is a Galois extension, let $G= \mathrm{Gal}(L|\fiel)$,
and suppose that the resulting group extension of $Q$ by $N$ 
is central.
The $Q$-action on $\mathrm U(K)$ passes to
an action of $Q$ on $\mathrm{coker}(\nu)$.
In terms of  classes of algebras of the kind $D(\tau,\eta)$,
the assignment to $D(\tau,\eta)$ of $D(\tau,{}^x\eta)$, as $x$ ranges over $Q$,
induces this action.

A little thought 
reveals that the following are equivalent: 
(i) {\em The algebra $D(\tau,\eta)$ 
 is $Q$-normal\/};
(ii) {\em the restriction 
$\Aut_L(D(\tau,\eta)) \to G$ is an epimorphism\/};
(iii) 
{\em
for $x \in Q=\mathrm{Gal}(K|\fiel)$, there is a unit 
$\vartheta \in \mathrm U(L)$ such that
$\nu(\vartheta) = {}^x \negthinspace\eta \eta^{-1}$. Furthermore, under the circumstances of {\rm (iii)},
the unit $\vartheta$ is unique up to multiplication 
by a unit in $K$.}

Condition (iii) plainly characterizes the members $[\eta]$ 
of the
subgroup $\mathrm{coker}(\nu)^Q$
of $Q$-invariants 
of the cokernel of the norm map $\nu \colon \mathrm U(L) \to \mathrm U(K)$.
Let $\sigma \colon Q \to \Aut(K)$ denote the Galois action.
The assignment to the class  $[\eta]\in \mathrm{coker}(\nu)^Q$ 
of the crossed $2$-fold extension $\mathrm e_{(D(\tau,\eta),\sigma)}$
of $Q$ by $\mathrm U(K)$
associated to a chosen  representative $\eta$
defines a map $t \colon \mathrm{coker}(\nu)^Q \to \Ho^3(Q,\mathrm U(K))$,
the {\em Teichmueller map\/} associated with the data.
Since for $\eta_1, \eta_2 \in \mathrm U(K)$, the tensor product algebra
$D(\tau,\eta_1) \otimes_KD(\tau,\eta_2)$
is the algebra of $(n \times n)$-matrices over $D(\tau,\eta_1 \eta_2)$,
the Teichmueller map $t$ is a homomorphism of abelian groups.

Up to this stage the discussion is elementary except, perhaps,
the quote of Hilbert's \lq\lq Satz 90\rq\rq. Now we
borrow some classical algebra:
Recall the  Brauer group $\mathrm B(K)$ of $K$ consists of
classes of central simple $K$-algebras, two such algebras being 
equivalent when they are matrix algebras over the same division algebra,
the inverse being induced by the assignment to an algebra of its 
opposite algebra.
A field $L|K$ {\em splits\/} the central simple $K$-algebra $A$
when $A \otimes _KL$ is $L$-isomorphic to a matrix algebra over $L$.
It is common to denote by $\mathrm B(L|K)$ the subgroup of Brauer classes
that are split by $L$.
For rings  more general than fields, the appropriate equivalence relation
is Morita equivalence.
The properties of being $Q$-normal and $Q$-equivariant
are properties of the Brauer classes,
the $Q$-action on $K$ induces an action of $Q$ on $\mathrm{Br}(K)$,
and the $Q$-invariants 
$\mathrm{Br}(K)^Q$ constitute the subgroup of Brauer classes of $Q$-normal
central simple $K$-algebras.
In terms of the canonical isomorphisms $\Ho^2(Q,\mathrm U(K)) \to
\mathrm{Br}(K|\fiel)$, $\Ho^2(G,\mathrm U(L)) \to
\mathrm{Br}(L|\fiel)$ and $\Ho^2(N,\mathrm U(L)) \to
\mathrm{Br}(L|K)$, with the notation $\mathrm{sc}$ for 
\lq scalar extension\rq, 
 the classical five term exact sequence in the cohomology of 
the (central) group extension of $Q$ by $N$ with coefficients in 
$\mathrm U(L)$ takes the form
\begin{equation}
\mathrm{Br}(K|\fiel) \rightarrowtail
\mathrm{Br}(L|\fiel) \stackrel{\mathrm{sc}}\rightarrow \mathrm{Br}(L|K)^Q\stackrel{t}
\rightarrow 
\Ho^3(Q, \mathrm U(K))
\stackrel{\mathrm{inf}}\rightarrow 
\Ho^3(G, \mathrm U(L)).
\label{fundfiv}
\end{equation}
To reconcile this sequence with the above remarks about 
the vanishing of the Teichmueller class of a
$Q$-equivariant
central $S$-algebra we note that, by \lq Galois descent\rq,
every $Q$-equivariant central simple $K$-algebra arises by scalar extension
from a central simple $\fiel$-algebra.

To arrive at explicit examples, let $K$ be
an algebraic number field
(a finite-dimensional extension of the field $\mathbb Q$ of rational numbers)
and,
as before, let $Q$ be a finite group of operators on $K$ and $\fiel  = K^Q$.
Let $m$ denote the l.c.m. of the local degrees
$[K_{\mathcat P}:\fiel_{\mathcat p}]$
as $\mathcat p$ ranges over the primes of 
$\fiel$ and
$\mathcat P$ over extensions thereof to $K$.
By
\cite[Theorem 3]{MR0025442}, the cokernel of 
scalar extension 
$\mathrm{sc}\colon\mathrm {Br}(\fiel) \to \mathrm {Br}(K)^Q$  
is a finite cyclic group of order $s= \tfrac {[K:\fiel]}m$.
Let $L=K(\zeta)$ be a  cyclotomic extension
having Galois group $N$ cyclic of order $n\geq 2$ (say), that is,
 $\zeta$ is a primitive $\ell$th root of unity for some $\ell$ 
prime to $n$, and
 the Galois group $N$ of order $n$ acts
faithfully and transitively on the 
 primitive $\ell$th roots of unity.
The field $L|\fiel$
coincides with  the composite field
$\fiel(\zeta)K$ in $L$,  and 
the canonical action of the pullback group
$ \mathrm{Gal}(\fiel(\zeta) |\fiel) 
\times_{\mathrm{Gal}(\fiel(\zeta)\cap K|\fiel)} Q$
on $\fiel(\zeta)K$ 
identifies this group with
a finite group of operators on $L=\fiel(\zeta)K$
having $\fiel$ as its fixed field. 
Hence
$L|\fiel$ is a Galois extension
having  Galois group $G$ canonially isomorphic
to the pullback group
$ \mathrm{Gal}(\fiel(\zeta) |\fiel) 
\times_{\mathrm{Gal}(\fiel(\zeta)\cap K|\fiel)} Q$, and
 $G$
 is a central extension of
$Q$ by the cyclic group 
$N=\mathrm{Gal}(K(\zeta)|K)
\cong\mathrm{Gal}(\fiel(\zeta)|(\fiel(\zeta)\cap K))$ of order $n$,
a split extension if and only if 
$\fiel(\zeta)\cap K = \fiel$.

Every class in $\mathrm{Br}(K)$ 
 has a cyclic cyclotomic  splitting field
but, beware, this only says that, for a central simple $K$-algebra $A$,
some cyclic cyclotomic  field 
splits a
matrix algebra over $A$.
On the other hand, it implies that
every class in $\mathrm{Br}(K)$ 
and in particular
in $\mathrm{Br}(K)^Q$ has a representative of the kind
$D(\tau,\eta)$.
Thus we may choose 
a member $\eta$ of $K$ 
and, for some $n\geq 2$, a cyclic degree $n$  cyclotomic Galois extension
$L = K(\zeta)$ of $K$
such that the
image $t [D(\tau,\eta]\in \Ho^3(Q, \mathrm U(K))$
of the  class $[D(\tau,\eta)] \in \mathrm{Br}(L|K)^Q$
of the $Q$-normal central $K$-algebra $D(\tau,\eta)$
generates the group $\Ho^3(Q, \mathrm U(K))$;
such a member $\eta$ is unique up to multiplication by some
$\lambda \in \fiel $ and by $\nu(\vartheta)$
for some $\vartheta \in L$.
From the exactness of \eqref{fundfiv} we deduce that
the Teichmueller map $t$  fits into the exact sequence
\begin{equation}
\mathrm{Br}(\fiel)\stackrel{\mathrm{sc}} 
\to \mathrm{Br}(K)^Q\stackrel{t}
\twoheadrightarrow 
\Ho^3(Q, \mathrm U(K))\cong \mathbb Z/s
\label{fit}
\end{equation}
but beware,
exactness only implies  that the class of a $Q$-normal algebra
having trivial Teichmueller class
arises by scalar extension, not necessarily the algebra itself.
An example of such an algebra that does not arise by scalar extension
while its class does is in 
 \cite{MR0002858}.
The images
$t [D(\tau,\eta)], t [D(\tau,\eta^2)], \ldots, t[D(\tau,\eta^{n-1})]
\in \Ho^3(Q, \mathrm U(K))$
exhaust the non-trivial members of the group $\Ho^3(Q, \mathrm U(K))$,
and the algebras
$D(\tau,\eta^j \lambda)$, for $1 \leq j \leq n-1$ and
$\lambda\in \fiel$, 
cover all Brauer classes of $Q$-normal algebras
split by $L$
with non-trivial Teichmueller class and, when we let $L$ vary,
we obtain  all Brauer classes of $Q$-normal algebras
with non-trivial Teichmueller class.
Thus, to get examples, all we need is
a Galois extension $K|\fiel$ having $s>1$.
While, in view of the Hilbert-Speiser Theorem,  
this is impossible when the Galois group $Q$ is cyclic,
for example,
the fields $K=\mathbb Q(\sqrt{13},\sqrt{17})$
or $K=\mathbb Q(\sqrt{2},\sqrt{17})$
have as Galois group  the four group and
$s=2$. 

Again we see that
 crossed modules having non-zero characteristic class 
and in particular non-extendible abstract kernels
abound.

\section{Outlook}
A topological group is a group in the category of topological spaces. 
Groups in the category $\mathrm{Cat}$ of small categories-equivalently,
categories internal to groups-, constitute a category, 
that of $2$-{\em groups\/}, and
there is
 an equivalence of categories
between crossed modules  and groups in 
$\mathrm{Cat}$,
observed by the Grothendieck school
in the mid 1960s (unpublished), see 
\cite[\S I.1.8 p.~29, \S 2.7 p.~58]{MR2841564}.
It is an interesting exercise to see how the Peiffer identities
fall out from this equivalence.

Crossed modules,
variants, and generalizations thereof are  nowadays
very lively in mathematics;
see  \cite{MR2841564},
\cite[Section 3]{MR4322147}
and the references there. 
The equivalence of crossed modules (in the category of groups) 
and $2$-groups
is relevant in string theory, see 
\cite[Section 3]{MR4322147} for references.
Suffice it to mention that when we pass from particles to strings, we add
an extra dimension, and replacing groups
by groups in $\mathrm{Cat}$
reflects this adding an extra dimension.

For a leisurely introduction and survey of the state of the art at the time
consider \cite{MR662431}. A particularly important work is
\cite{MR2841564}, with its special emphasis on foundational issues.

\section*{Acknowledgements}
I am indebted to AMS Notices editor Erica L. Flapan for encouraging me to
 compose this article
and to render it accessible to the non-expert,
and to the referee and to Jim Stasheff for a number of comments.
This work was supported in part by the Agence Nationale de la Recherche
under grant ANR-11-LABX-0007-01 (Labex CEMPI).

\bibliographystyle{alpha}

\def\cprime{$'$} \def\cprime{$'$} \def\cprime{$'$} \def\cprime{$'$}
  \def\cprime{$'$} \def\cprime{$'$} \def\cprime{$'$} \def\cprime{$'$}
  \def\dbar{\leavevmode\hbox to 0pt{\hskip.2ex \accent"16\hss}d}
  \def\cprime{$'$} \def\cprime{$'$} \def\cprime{$'$} \def\cprime{$'$}
  \def\cprime{$'$} \def\Dbar{\leavevmode\lower.6ex\hbox to 0pt{\hskip-.23ex
  \accent"16\hss}D} \def\cftil#1{\ifmmode\setbox7\hbox{$\accent"5E#1$}\else
  \setbox7\hbox{\accent"5E#1}\penalty 10000\relax\fi\raise 1\ht7
  \hbox{\lower1.15ex\hbox to 1\wd7{\hss\accent"7E\hss}}\penalty 10000
  \hskip-1\wd7\penalty 10000\box7}
  \def\cfudot#1{\ifmmode\setbox7\hbox{$\accent"5E#1$}\else
  \setbox7\hbox{\accent"5E#1}\penalty 10000\relax\fi\raise 1\ht7
  \hbox{\raise.1ex\hbox to 1\wd7{\hss.\hss}}\penalty 10000 \hskip-1\wd7\penalty
  10000\box7} \def\polhk#1{\setbox0=\hbox{#1}{\ooalign{\hidewidth
  \lower1.5ex\hbox{`}\hidewidth\crcr\unhbox0}}}
  \def\polhk#1{\setbox0=\hbox{#1}{\ooalign{\hidewidth
  \lower1.5ex\hbox{`}\hidewidth\crcr\unhbox0}}}
  \def\polhk#1{\setbox0=\hbox{#1}{\ooalign{\hidewidth
  \lower1.5ex\hbox{`}\hidewidth\crcr\unhbox0}}}

\end{document}